\documentclass{amsart}
\usepackage{amsfonts}

\usepackage{graphicx}
\usepackage{amscd}

\newtheorem{theorem}{Theorem}
\theoremstyle{plain}

\newtheorem{definition}{Definition}
\newtheorem{example}{Example}

\newtheorem{lemma}{Lemma}

\newtheorem{remark}{Remark}

\numberwithin{equation}{section}

\input{tcilatex}

\begin{document}
\title[Weighted kneading theory of unidimensional maps with holes]{Weighted kneading theory of unidimensional maps with holes}
\author{J. Leonel Rocha}
\address{Instituto Superior de Engenharia de Lisboa\\
Dep. de Eng. Qu\'{i}mica, Sec\c{c}\~{a}o de Matem\'{a}tica \\
Rua Conselheiro Em\'{i}dio Navarro \\
1949-014 Lisboa, Portugal}
\email{jrocha@math.ist.utl.pt}
\author{J. Sousa Ramos}
\curraddr{Instituto Superior T\'{e}cnico\\
Departamento de Matem\'{a}tica\\
Av. Rovisco Pais, 1 \\
1049-001 Lisboa, Portugal}
\email{sramos@math.ist.utl.pt}
\thanks{2000 \emph{Mathematics Subject Classification}: \emph{Primary}: 37E05,
37C30, 37B10; \emph{Secondary}: 37C45. \\
\emph{Keywords and phrases}: Kneading theory, Markov partitions, transfer
operator, topological pressure, topological entropy, Hausdorff dimension and
escape rate.}
\thanks{The first author is partially supported by PRODEP, Medida 5 - Ac\c{c}\~{a}o
5.3. Both thank FCT (Portugal) for having in part supported this work
through program POCTI. We are grateful to Nuno Martins for kindly reviewing
this paper and useful comments.}

\begin{abstract}
The purpose of this paper is to present a weighted kneading theory for
unidimensional maps with holes. We consider extensions of the kneading
theory of Milnor and Thurston to expanding discontinuous maps with holes and
introduce weights in the formal power series. This method allows us to
derive techniques to compute explicitly the topological entropy, the
Hausdorff dimension and the escape rate.
\end{abstract}

\maketitle

\section{Introduction}

Let $I\subset \mathbb{R}$ be a compact interval and $f=\left\{
f_{1},...,f_{n}\right\} $ be an iterated function system (IFS), a collection
of self-maps on $I$, defined by 
\begin{equation*}
f_{i}\left( x\right) :=\rho _{i}\,x+\varrho _{i},\,i=1,...,n
\end{equation*}
where for all $i$, $0<\left| \rho _{i}\right| <1$ and $\varrho _{i}\in 
\mathbb{R}$. Let $E$ be the corresponding self-similar set, the attractor.
We consider an IFS where the open set condition is satisfied, \cite{F}. If $%
f_{i}$ is monotone then is usual to see $E$ as the repeller of a linear
expanding map $F:\bigcup_{i=1}^{n}f_{i}\left( I\right) \rightarrow I$, which
will be denoted by $F=\left( F_{1},...,F_{n}\right) ,$ where 
\begin{equation*}
F_{i}\left( x\right) :=f_{i}^{-1}\left( x\right) \text{ if }x\in f_{i}\left(
I\right) .
\end{equation*}

We consider the piecewise linear map $F$ with a single hole, i.e., there is
an open subinterval $I_{h}\subset I$ with $I_{h}\neq \emptyset $ such that $%
I $ is the disjoint union of $I_{h}$ and $\bigcup_{i=1}^{n}\func{Im}\left(
f_{i}\right) $, (\cite{LS2} and \cite{LS3}). The points $x\in I_{h}$ will be
mapped out of $I$ and the same will happen to all the points $x\in
F^{-k}(I_{h})$ for $k\geq 1$. The set $\bigcup_{k}F^{-k}(I_{h})$ is open and
dense in $I$ and has full Lebesgue measure, \cite{BC}.

The hole and the set of $n$ laps of $F$ determines a partition $\mathcal{P}%
_{I}:=\left\{ I_{1},...,I_{h},...,I_{n}\right\} $ of the interval $I$.
Considering the orbits of the lateral limit points of the discontinuity
points and turning points, we define a Markov partition $\mathcal{P}%
_{I}^{\prime }$ of $I$.

The outline of the paper is as follows. In section 2, we develop a weighted
kneading theory to expanding discontinuous maps with holes. In this section
we give a brief presentation of the kneading theory associated to $F$. For
more details see \cite{MT} and for maps of the interval with holes see \cite
{LS1}. The new kneading approach using weights is inspired on \cite{LS3},
where we compute explicitly the escape rate that is characterized by a
conditional invariant measure. The weights introduced in the kneading theory
are defined by the inverses of the derivatives of the iterates of the
discontinuity points and turning points of $F$, associated to a real
parameter $\beta $. We consider the transfer or Perron-Frobenius operator $%
L_{\phi }$ associated to the map $F$ and to the Markov partition $\mathcal{P}%
_{I}^{\prime }$. The transfer operator has a matrix representation, which we
will denote by $Q_{\beta }$. This matrix can be viewed as the matrix of $%
L_{\phi }$ acting on a finite dimensional vector space of functions.

It is known that the spectrum of the transfer operator determines the
ergodicity of the dynamical system. In section 3, we will give an algorithm
to compute this spectrum and to relate the transfer operator with the
weighted kneading determinant, Theorem 2. To establish this relation we
introduce a weighted matrix\ $V_{\beta }$ and we use complexes and homology
with weights, Theorem 1. These results allow us to compute explicitly the
Hausdorff dimension, the escape rate and the topological entropy, Theorem 3.

We remark that, we can obtain the same results for a finite union of
disjoint holes $I_{h_{j}}\subset I$.

\section{Weighted kneading theory and subshifts of finite type}

Let $a_{i}$, with $i=1,...,n+1$, be the extreme points of the intervals in
the partition $\mathcal{P}_{I}$. These points correspond to the
discontinuity points and the turning points of the map $F$. We denote the
hole by $I_{h}=(a_{h},a_{h+1})$. Set 
\begin{equation*}
\begin{array}{l}
\left\{
a_{1}^{+},a_{2}^{-},a_{2}^{+},...,a_{h}^{-},a_{h}^{+},a_{h+1}^{-},a_{h+1}^{+},...,a_{n+1}^{-}\right\}
\\ 
\\ 
=\left\{ x^{\left( 1\right) },x^{\left( 2\right) },x^{\left( 3\right)
},...,x^{\left( 2h-2\right) },x^{\left( 2h-1\right) },x^{\left( 2\left(
h+1\right) -2\right) },x^{\left( 2\left( h+1\right) -1\right)
},...,x^{\left( 2n\right) }\right\} .
\end{array}
\end{equation*}
The orbit of each point $x^{\left( j\right) },$ with $j=1,\,2,...,\,2h-2,\,2%
\left( h+1\right) -1,$ $...,\,2n$, is defined by 
\begin{equation*}
o\left( x^{\left( j\right) }\right) :=\left\{ x_{k}^{\left( j\right)
}:x_{k}^{\left( j\right) }=F^{k}\left( x^{\left( j\right) }\right) ,\,k\in 
\mathbb{N}_{0}\right\} \text{.}
\end{equation*}
Concerning the itinerary of each point $x^{\left( j\right) }$ we will have 
\begin{equation*}
F^{k}\left( x^{\left( j\right) }\right) =x_{k}^{\left( j\right) }\text{ or }%
F^{k}\left( x^{\left( j\right) }\right) \in int\,(I_{h}).
\end{equation*}
In the first case, we have periodic, eventually periodic or aperiodic
orbits. While in the second, after a finite number of iterates, the
itinerary of the points lie in the hole.

To simplify the presentation, we consider the points $x^{\left( 1\right) }$
and $x^{\left( 2n\right) }$ as fixed points. To the orbit of each point $%
x^{\left( j\right) },$ with $j=2,...,\,2h-2,\,2\left( h+1\right)
-1,\,...,\,2n-1,$ we associate a sequence of symbols $S^{\left( j\right) }$
given by 
\begin{equation*}
S^{\left( j\right) }:=S_{0}^{\left( j\right) }S_{1}^{\left( j\right)
}...S_{k}^{\left( j\right) }...
\end{equation*}
where 
\begin{equation*}
S_{k}^{\left( j\right) }:=\left\{ 
\begin{array}{ll}
L & \text{if }F^{k}\left( x^{\left( j\right) }\right) \in I_{1} \\ 
&  \\ 
M_{s_{1}} & \text{if }F^{k}\left( x^{\left( j\right) }\right) \in
I_{s_{1}+1},\,s_{1}\in \left\{ 1,...,h-2\right\} \\ 
&  \\ 
H & \text{if }F^{k}\left( x^{\left( j\right) }\right) \in I_{h} \\ 
&  \\ 
M_{s_{2}} & \text{if }F^{k}\left( x^{\left( j\right) }\right) \in
I_{s_{2}+1},\text{\thinspace }s_{2}\in \left\{ h-1,...,n-2\right\} \\ 
&  \\ 
R & \text{if }F^{k}\left( x^{\left( j\right) }\right) \in I_{n}.
\end{array}
\right.
\end{equation*}
We denote by $\mathcal{A}$ the ordered set of $n+1$ symbols, corresponding
to the laps and the hole of $F$, i.e., 
\begin{equation*}
\mathcal{A}=\left\{
L,\,M_{1},...,\,M_{h-2},\,H,\,M_{h-1},...,\,M_{n-2},\,R\right\}
\end{equation*}
and according to the real line order 
\begin{equation*}
L\prec M_{1}\prec ...\prec M_{h-2}\prec H\prec M_{h-1}\prec ...\prec
M_{n-2}\prec R.
\end{equation*}
We designate by $\mathcal{A}^{\mathbb{N}}$ the space of all sequences of
symbols on the alphabet $\mathcal{A}$.

\begin{definition}
The kneading data for the map $F$ is the $(2n-4)$-tuple of symbolic
sequences 
\begin{equation*}
\left( S^{\left( 2\right) },...,S^{(2h-2)},S^{\left( 2\left( h+1\right)
-1\right) },...,S^{(2n-1)}\right) \in \mathcal{A}^{\mathbb{N}}\times 
\mathcal{A}^{\mathbb{N}}\times ...\times \mathcal{A}^{\mathbb{N}}.
\end{equation*}
\end{definition}

The kneading increments introduced in \cite{MT} are defined by formal power
series with coefficients in $\mathbb{Z}\left[ \left[ t\right] \right] $, the
subring of the ring $\mathbb{Q}\left[ \left[ t\right] \right] $. For maps of
the interval with holes and more than one discontinuity point and turning
points, we have several kneading increments, whose number depends on the
number of discontinuity points and turning points of the map $F$, \cite{LS1}%
. The kneading increments are defined by 
\begin{equation*}
\nu _{a_{i}}\left( t\right) :=\theta _{a_{i}^{+}}\left( t\right) -\theta
_{a_{i}^{-}}\left( t\right) .
\end{equation*}
In the case where $a_{i}$ is an extreme point of the hole $I_{h}$, the
increments are defined by 
\begin{equation*}
\nu _{a_{h}}\left( t\right) :=\theta _{a_{h}^{-}}\left( t\right) \text{ and }%
\nu _{a_{h+1}}\left( t\right) :=\theta _{a_{h+1}^{+}}\left( t\right)
\end{equation*}
where $\theta _{a_{i}}\left( t\right) $ is the invariant coordinate of each
symbolic sequence associated to the itinerary of each point $a_{i}$, with $%
1<i<n+1$. Each lateral invariant coordinate is defined by 
\begin{equation*}
\theta _{a_{i}^{\pm }}\left( t\right) :=\underset{x\rightarrow a_{i}^{\pm }}{%
\lim }\theta _{x}\left( t\right) =\underset{x\rightarrow a_{i}^{\pm }}{\lim }%
\sum_{k=0}^{\infty }\tau _{k}\,t^{k}\,S_{k}^{\left( j\right) }
\end{equation*}
where $\tau _{0}:=1$, $\tau _{k}:=\prod_{l=0}^{k-1}\varepsilon \left(
S_{l}^{\left( j\right) }\right) $, $k>0$ with 
\begin{equation*}
\varepsilon \left( S_{l}^{\left( j\right) }\right) :=\left\{ 
\begin{array}{ll}
\text{ \ }1 & \text{if }F^{\prime }\left( F^{l}\left( a_{i}^{\pm }\right)
\right) >0 \\ 
&  \\ 
\text{ \ }0 & \text{if }F^{\prime }\left( F^{l}\left( a_{i}^{\pm }\right)
\right) \in I_{h} \\ 
&  \\ 
-1 & \text{if }F^{\prime }\left( F^{l}\left( a_{i}^{\pm }\right) \right) <0
\end{array}
\right.
\end{equation*}
and $S_{k}^{\left( j\right) }$ is the sequence of symbols corresponding to
the orbits of $a_{i}^{\pm }$.

The increments $\nu _{a_{i}}\left( t\right) $,$\,$with $1<i<n+1$, can also
be written in the following way 
\begin{equation*}
\nu _{a_{i}}\left( t\right) =N_{i1}\left( t\right) \,L+N_{i2}\left( t\right)
\,M_{1}+...+N_{i(n-1)}\left( t\right) \,M_{n-2}+N_{in}\left( t\right) \,R
\end{equation*}
where the coefficients $N_{ij}\left( t\right) \in \mathbb{Z}\left[ \left[ t%
\right] \right] $ are the entries of a $n\times n$ matrix $N=\left[
N_{ij}\left( t\right) \right] $. This matrix is called the kneading matrix
associated to the map $F$. The kneading determinant is denoted by $D\left(
t\right) $.

Now we are going to present the main definition of this paper, that is the
characterization of the weighted invariant coordinates. This definition
allows us to construct a weighted kneading theory similar to the previous
one.

\begin{definition}
For the kneading data of the map $F$ the weighted invariant coordinate of
each point $a_{i}$, with $1<i<n+1$ and $\beta \in \mathbb{R}$, is defined by 
\begin{equation*}
\theta _{a_{i}^{\pm }}\left( t,\beta \right) :=\underset{x\rightarrow
a_{i}^{\pm }}{\lim }\theta _{x}\left( t,\beta \right) =\underset{%
x\rightarrow a_{i}^{\pm }}{\lim }\sum_{k=0}^{\infty }\tau _{k}\left(
a_{i}^{\pm }\right) \,t^{k}\,S_{k}^{\left( j\right) }
\end{equation*}
where $\tau _{0}\left( a_{i}^{\pm }\right) :=1$ and for $k>0$ 
\begin{equation*}
\tau _{k}\left( a_{i}^{\pm }\right) :=\prod_{l=0}^{k-1}\varepsilon \left(
S_{l}^{\left( j\right) }\right) \,\left| F^{\prime }\left( F^{l}\left(
a_{i}^{\pm }\right) \right) \right| ^{-\beta }.
\end{equation*}
\end{definition}

Note that, the derivative of the map $F$ satisfies the condition $\inf
\left| F_{i}^{\prime }\left( x\right) \right| >1$, on each interval $%
f_{i}\left( I\right) $. For each point $a_{i}$, with $%
i=1,...,h-1,h+2,...,n+1 $, the weighted kneading increment is defined by 
\begin{equation*}
\nu _{a_{i}}\left( t,\beta \right) :=\theta _{a_{i}^{+}}\left( t,\beta
\right) -\theta _{a_{i}^{-}}\left( t,\beta \right) .
\end{equation*}
For the extreme points of the hole $I_{h}$, the weighted increments are
defined by 
\begin{equation*}
\nu _{a_{h}}\left( t,\beta \right) :=\theta _{a_{h}^{-}}\left( t,\beta
\right) \text{ and }\nu _{a_{h+1}}\left( t,\beta \right) :=\theta
_{a_{h+1}^{+}}\left( t,\beta \right) .
\end{equation*}
Separating the terms associated to the symbols on the alphabet $\mathcal{A}$%
, the weighted increments $\nu _{a_{i}}\left( t,\beta \right) $ are written
in the following way 
\begin{equation*}
\nu _{a_{i}}\left( t,\beta \right) =N_{i1}\left( t,\beta \right)
\,L+N_{i2}\left( t,\beta \right) \,M_{1}+...+N_{i(n-1)}\left( t,\beta
\right) \,M_{n-2}+N_{in}\left( t,\beta \right) \,R.
\end{equation*}
The coefficients $N_{ij}\left( t,\beta \right) \in \mathbb{R}\left[ \left[
t,\beta \right] \right] $ are the entries of a $n\times n$ matrix 
\begin{equation*}
N\left( t,\beta \right) :=\left[ N_{ij}\left( t,\beta \right) \right]
\end{equation*}
which we will call the weighted kneading matrix associated to $F$. The
determinant of this matrix will be called the weighted kneading determinant
and will be denoted by $D\left( t,\beta \right) $.

\begin{remark}
To an eventually periodic orbit of a point $x^{\left( j\right) }$
represented by 
\begin{equation*}
S_{0}^{\left( j\right) }...S_{p-1}^{\left( j\right) }\left( S_{p}^{\left(
j\right) }...S_{p+q-1}^{\left( j\right) }\right) ^{\infty }
\end{equation*}
it corresponds the weighted cyclotomic polynomial 
\begin{equation*}
1-\prod_{l=0}^{q-1}\varepsilon \left( S_{p+l}^{\left( j\right) }\right)
\,\left| F^{\prime }\left( F^{p+l}\left( x^{\left( j\right) }\right) \right)
\right| ^{-\beta }\,t^{q}
\end{equation*}
where $q$ is the period of the orbit. If the orbit is periodic then the
weighted cyclotomic polynomial is 
\begin{equation*}
1-\tau _{q}\left( x^{\left( j\right) }\right) \,t^{q}.
\end{equation*}
\end{remark}

Now, let 
\begin{equation*}
\left\{ b_{1},...,b_{m+1}\right\} :=\left\{ o\left( x^{\left( j\right)
}\right) :j=1,...,\,2h-2,\,2\left( h+1\right) -1,\,...,\,2n\right\}
\end{equation*}
be the set of the points correspondent to the orbits of the lateral limit
points of the discontinuity points and turning points, ordered on the
interval $I$. This set allows us to define a subpartition $\mathcal{P}%
_{I}^{\prime }$ of $\mathcal{P}_{I}=\left\{
I_{1},...,I_{h},...,I_{n}\right\} .$ The subpartition 
\begin{equation*}
\mathcal{P}_{I}^{\prime }:=\left\{ J_{1},...,J_{m}\right\}
\end{equation*}
with $m\geq n$ determines a Markov partition of the interval $I$. Note that
the hole is an element of the Markov partition. Note also that $F$
determines $\mathcal{P}_{I}^{\prime }$ uniquely, but the converse is not
true.

The IFS $f$ induces a subshift of finite type whose $m\times m$ transition
matrix $A=[a_{ij}]$ is defined by 
\begin{equation*}
a_{ij}:=\left\{ 
\begin{array}{ll}
1 & \text{if }F\left( int\,J_{j}\right) \supseteq int\,J_{i} \\ 
0 & \text{otherwise.}
\end{array}
\right.
\end{equation*}
We remark that if there exists $k$ points $b_{i}$ such that $b_{i}\in
int\,I_{h}$, with $1<i<m+1$, then the matrix $A$ has $k+1$ columns with all
elements equal to zero, correspondent to the hole.

We denote this subshift by $\left( \Sigma _{A},\sigma \right) $, where $%
\sigma $ is the shift map on $\Sigma _{m}^{\mathbb{N}}$ defined by $\sigma
\left( x_{1}x_{2}...\right) :=x_{2}x_{3}...$, with $\Sigma _{m}:=\left\{
1,...,m\right\} $ correspondent to the $m$ states of the subshift.
Concerning this subshift $\left( \Sigma _{A},\sigma \right) $ and the
associated Markov partition $\mathcal{P}_{I}^{\prime }$, we consider a
Lipschitz function $\phi :I\rightarrow \mathbb{R}$, defined by 
\begin{equation*}
\phi :=\left\{ \phi _{i}:J_{i}\rightarrow \mathbb{R},\,1\leq i\leq m\right\}
\end{equation*}
where 
\begin{equation*}
\phi _{i}\left( x\right) :=-\beta \text{\thinspace }\varphi _{i}\left(
x\right) \text{ and }\varphi _{i}\left( x\right) :=\log \left| F_{i}^{\prime
}\left( x\right) \right| \text{, with }\beta \in \mathbb{R}.
\end{equation*}
This function is a weight for the dynamical system associated to the
subshift, depending on the real parameter $\beta $ (compare with \cite{R1}).

Let $\mathcal{L}^{1}\left( I\right) $ be the set of all Lebesgue integrable
functions on $I$. The transfer operator $L_{\phi }:\mathcal{L}^{1}\left(
I\right) \rightarrow \mathcal{L}^{1}\left( I\right) $, associated with $F$
and $\mathcal{P}_{I}$, is defined by 
\begin{equation*}
\begin{array}{lll}
\left( L_{\phi }\,g\right) \left( x\right) & := & \sum\limits_{j=1}^{n}\exp
\phi \left( F_{j}^{-1}\left( x\right) \right) \text{\thinspace }g\left(
F_{j}^{-1}\left( x\right) \right) \,\chi _{F\left( int\,I_{j}\right) } \\ 
&  &  \\ 
& = & \sum\limits_{y:Fy=x}\exp \phi \left( y\right) \text{\thinspace }%
g\left( y\right)
\end{array}
\end{equation*}
where $\chi _{I_{j}}$ is the characteristic function of $I_{j}$. Note that,
by definition of $F$, $F_{j}^{-1}\left( x\right) =f_{j}\left( x\right) $
with $x\in F\left( I_{j}\right) $. Note also that, for any Borel subset $%
J\subset I,$ we have 
\begin{equation*}
F^{-1}\left( J\right) =\bigcup\limits_{j=1}^{n}f_{j}\left( F\left(
I_{j}\right) \cap J\right)
\end{equation*}
where the sets $\left\{ f_{j}\left( F\left( I_{j}\right) \cap J\right)
\right\} _{j=1}^{n}$ are mutually disjoint. Depending on $J$, the set $%
f_{j}\left( F\left( I_{j}\right) \cap J\right) $ can be empty.

Now, we will restrict our attention to the transfer operator associated with 
$F$ and to the Markov partition $\mathcal{P}_{I}^{\prime }$. Given $J_{i}\in 
\mathcal{P}_{I}^{\prime }$, let $Y_{1i},...,Y_{ki}$ be the preimages of $%
J_{i}$ under $F$, i.e., 
\begin{equation*}
F_{j}\left( int\,Y_{ji}\right) =int\,J_{i}\text{ with }1\leq j\leq k\text{
and }k\leq m.
\end{equation*}
Then, we can define continuous maps $f_{j}\left| _{J_{i}}\right. :=\Psi
_{ji}:J_{i}\rightarrow Y_{ji}$, that correspond to the IFS $f$ restricted to
the interval $J_{i}$, such that $y_{j}=\Psi _{ji}\left( x\right) $ are the
preimages of $x\in J_{i}$. Thus, for each $x\in J_{i}$ we have 
\begin{equation}
\left( L_{\phi _{j}}\,g\right) \left( x\right) =\sum\limits_{j=1}^{k}\exp
\phi _{j}\left( \Psi _{ji}\left( x\right) \right) \text{\thinspace }g\left(
\Psi _{ji}\left( x\right) \right) \,\delta \left( \Psi _{ji}\left( x\right)
\right)  \label{FK1.1}
\end{equation}
where 
\begin{equation*}
\delta \left( \Psi _{ji}\left( x\right) \right) :=\left\{ 
\begin{array}{ll}
1 & \text{if }\Psi _{ji}\left( x\right) =y_{j} \\ 
0 & \text{otherwise.}
\end{array}
\right.
\end{equation*}
Nevertheless, for each interval $J_{j}\in \mathcal{P}_{I}^{\prime }$ we
consider 
\begin{equation*}
f_{j}\left( x\right) :=F^{-1}\left| _{J_{j}}\right. \left( x\right) \text{
for }x\in F\left( int\,J_{j}\right) =\bigcup\limits_{a_{ij}\neq 0}int\,J_{i}
\end{equation*}
where $a_{ij}$ are the entries of the transition matrix $A,$ with $1\leq
i,\,j\leq m.$ By formula (\ref{FK1.1}), we can write 
\begin{equation}
\left( L_{\phi _{j}}\,g\right) \left( x\right) =\sum\limits_{j=1}^{m}\exp
\phi _{j}\left( f_{j}\left( x\right) \right) \text{\thinspace }g\left(
f_{j}\left( x\right) \right) \,\chi _{F\left( int\,J_{j}\right) }.
\label{FK1.2}
\end{equation}

In this paper we consider a class of one-dimensional transformations that
are piecewise linear Markov transformations. Consequently, the transfer
operator has the following matrix representation. Let $\mathcal{C}$ be the
class of all functions that are piecewise constant on the partition $%
\mathcal{P}_{I}^{\prime }$. Thus, 
\begin{equation*}
g\in \mathcal{C}\text{ iff }g=\sum\limits_{j=1}^{m}\pi _{j}\,\chi _{J_{j}}
\end{equation*}
for some constants $\pi _{1},...,\pi _{m}$. We remark that $g$ will also be
represented by the column vector $\pi _{g}=\left( \pi _{1},...,\pi
_{m}\right) ^{T}$. Using the formula (\ref{FK1.2}) and considering $g\in 
\mathcal{C}$ with $g=\chi _{J_{k}}$, for some $1\leq k\leq m$, the transfer
operator $L_{\phi }$ has the following matrix characterization 
\begin{equation*}
L_{\phi }\,g=Q_{\beta }\,\pi _{g}
\end{equation*}
for the weighted dynamical system associated to $\left( \Sigma _{A},\sigma
\right) $. If $D_{\beta }$ is the diagonal matrix defined by 
\begin{equation*}
D_{\beta }:=\left( \exp \varphi _{1}^{-\beta },...,\exp \varphi _{m}^{-\beta
}\right)
\end{equation*}
and $A$ is the transition matrix, then the matrix $Q_{\beta }$ is the $%
m\times m$ weighted transition matrix defined by 
\begin{equation*}
Q_{\beta }:=A\,D_{\beta }=[q_{ij}].
\end{equation*}
The entries of this matrix are 
\begin{equation*}
q_{ij}:=\frac{a_{ij}}{\left| F_{j}^{\prime }\right| ^{\beta }}
\end{equation*}
where the derivative $F_{j}^{\prime }$ is evaluated on the interval $J_{j}$
of the partition $\mathcal{P}_{I}^{\prime }$. We refer to \cite{R2}, \cite
{BG} and the references therein to other important spectral properties of
the transfer operator, and \cite{F} for this operator with respect to the
cookie-cutter system. In \cite{LS3}, we use the matrix $Q_{\beta }$ with $%
\beta =1$ to compute the escape rate and the conditional invariant measure
which generates the unique invariant probability measure.

There is an isomorphism between $\left( \Sigma _{A},\sigma \right) $ and $%
\left( \mathcal{P}_{I}^{\prime },F\right) $, \cite{W}. If $w=\left(
i_{0}i_{1}...\right) $ and $w^{\prime }=\left( i_{0}^{\prime }i_{1}^{\prime
}...\right) $ are two points of $\Sigma _{A}$, then we consider the Markov
partition defined by $\mathcal{P}_{I}^{\prime }:=\Sigma _{A}/\sim $ where $%
w\sim w^{\prime }$ if and only if $i_{0}=i_{0}^{\prime }$. Using this
isomorphism, we consider the trace of the transfer operator defined by 
\begin{equation*}
Tr\,L_{\phi }:=\sum\limits_{x\in \text{Fix}(F)}\exp \phi \left( x\right)
\end{equation*}
where Fix$(F)$ denotes the set of fixed points of $F.$ We consider the
pressure function of $\phi \left( x\right) =\log \left| F^{\prime }\left(
x\right) \right| ^{-\beta }$ as $\beta $ varies, $P\left( \beta \right) $,
defined by 
\begin{equation}
\begin{array}{lll}
P\left( \beta \right) & := & \underset{k\rightarrow \infty }{\lim }\frac{1}{k%
}\log \sum\limits_{x\in \text{Fix}(F^{k})}\left| \left( F^{k}\right)
^{\prime }\left( x\right) \right| ^{-\beta } \\ 
&  &  \\ 
& = & \underset{k\rightarrow \infty }{\lim }\frac{1}{k}\log \left(
Tr\,Q_{\beta }^{k}\right) \\ 
&  &  \\ 
& = & \log \left( \lambda _{\beta }\right)
\end{array}
\label{FK1.4}
\end{equation}
where Fix$(F^{k})$ denotes the set of fixed points of $F^{k}$, (\cite{F}, 
\cite{W} and \cite{R3}). Thus, $\exp P\left( \beta \right) $ is the largest
eigenvalue $\lambda _{\beta }$ of the transfer operator $L_{\phi }$, which
is equal to the spectral radius of the matrix $Q_{\beta }$, \cite{R2}.
Nevertheless, it is known that the pressure can be characterized by the
variational principle 
\begin{equation*}
P\left( \beta \right) =h_{\mu _{\beta }}\left( F\right) -\beta \,\chi _{\mu
_{\beta }}\left( F\right)
\end{equation*}
as the supremum over all invariant probability measures on $E$. In this
case, the supremum is attained by the weighted Markov measure $\mu _{\beta }$%
, i.e., the measure $\mu _{\beta }$ is the unique measure that maximizes
this expression. See \cite{LS3} for the definition of $\mu _{\beta }$, the
weighted metric entropy $h_{\mu _{\beta }}\left( F\right) $ and the weighted
Lyapunov exponent $\chi _{\mu _{\beta }}\left( F\right) $ with respect to
this measure.

We remark that the weighted zeta function for a weighted subshift of finite
type is given by 
\begin{equation*}
\zeta \left( t,\beta \right) =\frac{1}{\det \left( I-t\,Q_{\beta }\right) }.
\end{equation*}
For more discussions about the zeta function for a subshift of finite type
without weights see \cite{BL}, and for another approach with weights see 
\cite{M} and \cite{BR}.

\section{Complexes and homology with weights}

As above described, on the set $\left\{ b_{1},...,\,b_{m+1}\right\} $ there
exist $2n-2$ points 
\begin{equation}
\left\{ x^{\left( 1\right) },...,x^{\left( 2h-2\right) },x^{\left(
2h+1\right) },...,x^{\left( 2n\right) }\right\} .  \label{FK2.0}
\end{equation}
Concerning the set of points $\left\{ a_{1},...,\,a_{n+1}\right\} $, we
consider the respective $2n-2$ lateral points, as in section 2. Let $p$ be
the number of points $b_{i}$ outside of (\ref{FK2.0}). We denote this set by 
$\mathcal{G}$. Set $2n+p-2=q$ and 
\begin{equation*}
\begin{array}{l}
\left\{
a_{1}^{+},\,a_{2}^{-},\,a_{2}^{+},...,\,a_{h}^{-},\,a_{h}^{+},\,a_{h+1}^{-},%
\,a_{h+1}^{+},...,\,a_{n+1}^{-}\right\} \cup \mathcal{G} \\ 
\\ 
=\left\{ o\left( x^{\left( j\right) }\right) :j=1,...,\,2h-2,\,2\left(
h+1\right) -1,...,\,2n\right\} \cup \left\{ x^{\left( 2h-1\right)
},\,x^{\left( 2(h+1)-2\right) }\right\} \\ 
\\ 
=\left\{ x^{\left( 1\right) },\,x_{1}^{\left( 1\right)
},...,\,x_{k_{1}}^{\left( 1\right) },...,\,x^{\left( 2n\right)
},\,x_{1}^{\left( 2n\right) },...,\,x_{k_{2n}}^{\left( 2n\right) }\right\}
\cup \left\{ x^{\left( 2h-1\right) },\,x^{\left( 2(h+1)-2\right) }\right\}
\\ 
\\ 
=\left\{ z^{\left( 1\right) },\,z^{\left( 2\right) },...,\,z^{\left(
q-2\right) },\,x^{\left( 2h-1\right) },\,x^{\left( 2(h+1)-2\right) }\right\}
\\ 
\\ 
=\left\{ z^{\left( 1\right) },\,z^{\left( 2\right) },...,\,z^{\left(
q\right) }\right\} .
\end{array}
\end{equation*}
We denote by $\left\{ y^{\left( 1\right) },\,y^{\left( 2\right)
},...,\,y^{\left( q\right) }\right\} $ the above set of points now ordered
on the interval $I$. On the set $\left\{ y^{\left( 1\right) },\,y^{\left(
2\right) },...,\,y^{\left( q\right) }\right\} $ there exist pairs of
consecutive points $y^{\left( k\right) },y^{\left( k+1\right) }$, with $%
2\leq k\leq q-2$, correspondent to some $a_{i}^{\pm }$ with $i=2,...,n$.
According to the above conditions, we define a permutation $\rho $ by 
\begin{equation*}
\left( z^{\left( 1\right) },\,z^{\left( 2\right) },...,\,z^{\left( q\right)
}\right) \rightarrow \left( z^{\rho \left( 1\right) },\,z^{\rho \left(
2\right) },...,\,z^{\rho \left( q\right) }\right) =\left( y^{\left( 1\right)
},\,y^{\left( 2\right) },...,\,y^{\left( q\right) }\right) .
\end{equation*}

Let $C_{0}$ be the vector space of $0$-chains spanned by the points $%
y^{\left( 1\right) },...,y^{\left( q\right) }$ and $C_{1}$ be the vector
space of $1$-chains spanned by the intervals of the partition $\mathcal{P}%
_{I}^{\prime }$. The border map $\partial :C_{1}\rightarrow C_{0}$ is
defined by 
\begin{equation*}
\partial \left( J_{j}\right) :=y^{\left( i+1\right) }-y^{\left( i\right) }
\end{equation*}
with $1\leq j\leq m$ and $1\leq i\leq q-1.$ We designate by $B=\left[ b_{ij}%
\right] $ the $q\times m$ incidence matrix of the graph defined by $\left(
C_{0},C_{1},\partial \right) $, i.e., if $y^{\left( i\right) }$ is the lower
extreme point and $y^{\left( i+1\right) }$ is the upper extreme point of the
interval $J_{j}$ then $b_{ij}:=-1$ and $b_{i+1,j}:=1$, the remaining entries
are zero.

Let $\eta :C_{0}\rightarrow C_{0}$ be the map that describes the transition
between the points $y^{\left( 1\right) },y^{\left( 2\right) },...,y^{\left(
q\right) }$ and checks the existence of turning points and discontinuity
points between $y^{\left( j\right) }$ and $F\left( y^{\left( j\right)
}\right) $, with $1\leq j\leq q$. This map is represented by the following
weighted matrix $V_{\beta }$.

\begin{definition}
The $q\times q$ weighted matrix $V_{\beta }=\left[ v_{ij}\right] $,
associated to the map $F$, is defined by 
\begin{equation*}
\begin{array}{ll}
v_{ij}:=\varepsilon \left( y^{\left( j\right) }\right) \,\left| F^{\prime
}\left( y^{\left( j\right) }\right) \right| ^{-\beta } & \text{if }F\left(
y^{\left( j\right) }\right) =y^{\left( i\right) }\text{, where }\varepsilon
\left( y^{\left( j\right) }\right) =sign\left( F^{\prime }\left( y^{\left(
j\right) }\right) \right) \\ 
&  \\ 
v_{kj}:=v_{ij},\,v_{k+1,j}:=-v_{ij} & \text{if }y^{\left( i\right)
}>y^{\left( j\right) }\text{ and }j\leq k<i \\ 
&  \\ 
v_{k-1,j}:=-v_{ij},\,v_{k,j}:=v_{ij} & \text{if }y^{\left( i\right)
}<y^{\left( j\right) }\text{ and }i<k\leq j
\end{array}
\end{equation*}
where the pairs of consecutive points $y^{\left( k\right) },\,y^{\left(
k+1\right) }$ with $j\leq k<i$ (resp. $y^{\left( k-1\right) },\,y^{\left(
k\right) }$ with $i<k\leq j$) are associated to the turning points and
discontinuity points of $F$ with $y^{\left( j\right) }\leq y^{\left(
k\right) },\,y^{\left( k+1\right) }\leq y^{\left( i\right) }$ (resp. $%
y^{\left( i\right) }\leq y^{\left( k-1\right) },\,y^{\left( k\right) }\leq
y^{\left( j\right) }$). All the remaining entries in $V_{\beta }$ are zero.
\end{definition}

We consider that the points $y^{\left( i\right) }$, with $1\leq i\leq q$,
are represented by $y^{\left( i\right) }=\left( 0,...,0,1,0,...,0\right)
^{T} $, where $1$ is in the $i^{\text{th}}$-position. The above weighted
matrices are related by the next result.

\begin{lemma}
The next diagram is commutative 
\begin{equation*}
\begin{array}{ccccc}
&  & B &  &  \\ 
& C_{1} & \longrightarrow & C_{0} &  \\ 
Q_{\beta } & \downarrow &  & \downarrow & V_{\beta } \\ 
& C_{1} & \longrightarrow & C_{0} &  \\ 
&  & B &  & 
\end{array}
.
\end{equation*}
\end{lemma}

\begin{proof}
Let the intervals $J_{j}$, with $1\leq j\leq m$, be represented by the
column vector $J_{j}=\left( 0,...,0,\chi _{J_{j}},0,...,0\right) ^{T}$,
where $\chi _{J_{j}}$ is in the $j^{\text{th}}$-position, i.e., as a
function in $\mathcal{C}$. Then, one has 
\begin{equation*}
\begin{array}{lll}
\left( B\,Q_{\beta }\right) \,\left( J_{j}\right) & = & \left(
B\,A\,D_{\beta }\right) \,\left( J_{j}\right) \\ 
&  &  \\ 
& = & \left( B\,A\right) \,\left( e^{\varphi _{j}^{-\beta }}\,J_{j}\right)
\\ 
&  &  \\ 
& = & B\,\left( \sum\limits_{a_{kj}\neq 0}\left| F^{\prime }\left(
J_{j}\right) \right| ^{-\beta }\,J_{k}\right) \\ 
&  &  \\ 
& = & \sum\limits_{a_{kj}\neq 0}\left| F^{\prime }\left( J_{j}\right)
\right| ^{-\beta }\,B\,\left( J_{k}\right) \\ 
&  &  \\ 
& = & \sum\limits_{u_{lj}\neq 0}\left| F^{\prime }\left( J_{j}\right)
\right| ^{-\beta }\,\left( y^{\left( l+1\right) }-y^{\left( l\right) }\right)
\end{array}
\end{equation*}
where $u_{lj}$ are the nonzero elements of the $j^{\text{th}}$ column of the
matrix $BA$. The above equalities make a description of the transition of
the interval $J_{j}$ by the border of the intervals $J_{k}$ such that $%
F\left( int\,J_{j}\right) \supseteq int\,J_{k}$, weighted by $\left|
F^{\prime }\left( J_{j}\right) \right| ^{-\beta }$.

On the other hand, we have 
\begin{equation*}
V_{\beta }\,B\left( J_{j}\right) =V_{\beta }\,\left( y^{\left( i+1\right)
}-y^{\left( i\right) }\right) =V_{\beta }\,\left( y^{\left( i+1\right)
}\right) -V_{\beta }\,\left( y^{\left( i\right) }\right) .
\end{equation*}
Consider that $y^{\left( i\right) }$ is a point associated to a turning
point or to a discontinuity point of $F$ and $F\left( y^{\left( i\right)
}\right) =y^{\left( s\right) }$, with $s<i$. Consequently, 
\begin{equation*}
v_{si}=\varepsilon \left( y^{\left( i\right) }\right) \,\left| F^{\prime
}\left( y^{\left( i\right) }\right) \right| ^{-\beta },\,v_{i-1,i}=-v_{si}%
\text{ and }v_{ii}=v_{si}.
\end{equation*}
If there exist $z_{1}$ turning points or discontinuity points between $%
y^{\left( s\right) }$ and $y^{\left( i-1\right) }$, then we have pairs of
consecutive points $y^{\left( k_{l}-1\right) },y^{\left( k_{l}\right) }$,
with $s<k_{l}-1,k_{l}<i-1$ and $1\leq l\leq z_{1}$, such that 
\begin{equation*}
v_{k_{l}-1,i}=-v_{si}\text{ and }v_{k_{l}i}=v_{si}.
\end{equation*}
Suppose that $y^{\left( i+1\right) }\in int\,I_{p}$ with $1\leq p\leq n$,
i.e., $y^{\left( i+1\right) }\in \mathcal{G}$ and $F\left( y^{\left(
i+1\right) }\right) =y^{\left( r\right) }$, with $r>i+1.$ In this case, we
have 
\begin{equation*}
v_{i+1,i+1}=0\text{ and }v_{r,i+1}=\varepsilon \left( y^{\left( i+1\right)
}\right) \,\left| F^{\prime }\left( y^{\left( i+1\right) }\right) \right|
^{-\beta }.
\end{equation*}
Similarly, if there exists $z_{2}$ turning points or discontinuity points
between $y^{\left( i+1\right) }$ and $y^{\left( r\right) }$, then we have
pairs of points $y^{\left( k_{w}\right) },\,y^{\left( k_{w}+1\right) }$,
with $i+1<k_{w},\,k_{w}+1<r$ and $1\leq w\leq z_{2}$ such that 
\begin{equation*}
v_{k_{w}i}=v_{r,i+1}\text{ and }v_{k_{w}+1,i}=-v_{r,i+1}.
\end{equation*}
As the weight is constant on each interval $J_{j}$, we get 
\begin{equation*}
V_{\beta }\,\left( y^{\left( i+1\right) }\right) -V_{\beta }\,\left(
y^{\left( i\right) }\right) =\left| F^{\prime }\left( J_{j}\right) \right|
^{-\beta }\,\sum\limits_{l=s}^{r-1}\left( y^{\left( l+1\right) }-y^{\left(
l\right) }\right)
\end{equation*}
where the pairs of points $y^{\left( l\right) },\,y^{\left( l+1\right) }$
lie in the set 
\begin{equation*}
\left\{ y^{\left( s\right) },y^{\left( s+1\right) },...,y^{\left(
k_{l}\right) },y^{\left( k_{l}+1\right) },...,y^{\left( i-1\right)
},y^{\left( i\right) },...,y^{\left( k_{w}\right) },y^{\left( k_{w}+1\right)
},...,y^{\left( r-1\right) },y^{\left( r\right) }\right\}
\end{equation*}
and describe the border of the intervals $J_{k}$ such that $F\left(
int\,J_{j}\right) \supseteq int\,J_{k}$. From this it follows that 
\begin{equation*}
\left| F^{\prime }\left( J_{j}\right) \right| ^{-\beta
}\,\sum\limits_{l=s}^{r-1}\left( y^{\left( l+1\right) }-y^{\left( l\right)
}\right) =\sum\limits_{a_{kj}\neq 0}\left| F^{\prime }\left( J_{j}\right)
\right| ^{-\beta }\,B\left( J_{k}\right) .
\end{equation*}
The proof of the remaining cases is similar, according to the behavior of $F$
and the above definition.
\end{proof}

Let $H_{0}:=C_{0}/B_{0}$, where $B_{0}=\partial \left( C_{1}\right) $ is a
subspace of $C_{0}$. Note that two consecutive laps without a discontinuity
point between them are considered as two connected components. The map $%
\zeta :C_{0}\rightarrow H_{0}$ associates to each point $y^{\left( i\right)
} $, with $1\leq i\leq q$, the respective interval $I_{j}$, with $1\leq
j\leq n $. This map is represented by the $q\times n$ matrix $U=\left[ u_{ij}%
\right] $, where $u_{ij}:=1$ if the point $y^{\left( i\right) }$ lies in $%
I_{j}$ and all the remaining entries of the matrix are zero.

\begin{lemma}
If $y^{\left( i_{1}\right) }$ and $y^{\left( i_{2}\right) }$ are two points
on the interval $I_{j}$ with $1\leq j\leq n$, then we have 
\begin{equation*}
U\,V_{\beta }\left( y^{\left( i_{1}\right) }\right) ^{T}=U\,V_{\beta }\left(
y^{\left( i_{2}\right) }\right) ^{T}.
\end{equation*}
\end{lemma}

\begin{proof}
Let $y^{\left( i_{1}\right) }$ and $y^{\left( i_{2}\right) }$ be two points
on $I_{j}$. Consider that, 
\begin{equation*}
B\,\left( J_{k}\right) =y^{\left( i_{2}\right) }-y^{\left( i_{1}\right)
},\,int\,J_{k}\subseteq int\,I_{j},\,F\left( y^{\left( i_{1}\right) }\right)
\in I_{j^{\prime }}\text{ and }F\left( y^{\left( i_{2}\right) }\right) \in
I_{j^{\prime \prime }}
\end{equation*}
with $1\leq k\leq m$ and $j^{\prime }<j<j^{\prime \prime }$. If the points $%
y^{\left( i_{1}\right) }$ and $y^{\left( i_{2}\right) }$ have the same
behavior under $F$ as $y^{\left( i\right) }$ and $y^{\left( i+1\right) }$ in
the proof of the above proposition, respectively, then 
\begin{equation*}
\begin{array}{lll}
U\,V_{\beta }\,\left( y^{\left( i_{1}\right) }\right) ^{T} & = & U\,\left(
0,...,0,v_{si_{1}},...,v_{k_{l}-1,i_{1}},v_{k_{l}i_{1}},...,v_{i_{1}-1,i_{1}},v_{i_{1}i_{1}},0,...,0\right) ^{T}
\\ 
&  &  \\ 
& = & U\,\varepsilon \left( y^{\left( i_{1}\right) }\right) \,\left|
F^{\prime }\left( y^{\left( i_{1}\right) }\right) \right| ^{-\beta }\left(
0,...,0,1,...,-1,1,...,-1,1,0,...,0\right) ^{T}.
\end{array}
\end{equation*}
Note that there exist $z_{1}$ pairs of consecutive points $y^{\left(
k_{l}-1\right) },\,y^{\left( k_{l}\right) }$ between $y^{\left( s\right) }$
and $y^{\left( i_{1}-1\right) },$ with $1\leq l\leq z_{1}.$ These points
define the border of the $p$ intervals $I_{j_{t}}$, with $j^{\prime
}<j_{t}<j $ and $1\leq t\leq p$. This implies that 
\begin{equation*}
\begin{array}{lll}
B\,\left( I_{j^{\prime }}\right) & = & y^{\left( k_{1}-1\right) }-y^{\left(
s\right) } \\ 
&  &  \\ 
B\,\left( I_{j_{1}}\right) & = & y^{\left( k_{2}-1\right) }-y^{\left(
k_{1}\right) } \\ 
\text{ \ \ \ }\vdots &  &  \\ 
B\,\left( I_{j}\right) & = & y^{\left( i_{1}\right) }-y^{\left(
i_{1}-1\right) }.
\end{array}
\end{equation*}
Thus, by definition of $U$ we get 
\begin{equation*}
U\,V_{\beta }\,\left( y^{\left( i_{1}\right) }\right) ^{T}=\left(
0,...,0,\varepsilon \left( y^{\left( i_{1}\right) }\right) \,\left|
F^{\prime }\left( y^{\left( i_{1}\right) }\right) \right| ^{-\beta
},0,...,0\right) ^{T}
\end{equation*}
where $\varepsilon \left( y^{\left( i_{1}\right) }\right) \,\left| F^{\prime
}\left( y^{\left( i_{1}\right) }\right) \right| ^{-\beta }$ is in the $j^{%
\text{th}}$-position.

On the other hand, considering the point $y^{\left( i_{2}\right) }$ we can
write 
\begin{equation*}
\begin{array}{lll}
U\,V_{\beta }\,\left( y^{\left( i_{2}\right) }\right) ^{T} & = & U\,\left(
0,\ldots ,0,v_{k_{w}i_{2}},v_{k_{w}+1,i_{2}},...,v_{ri_{2}},0,...,0\right)
^{T} \\ 
&  &  \\ 
& = & U\,\varepsilon \left( y^{\left( i_{2}\right) }\right) \,\left|
F^{\prime }\left( y^{\left( i_{2}\right) }\right) \right| ^{-\beta }\,\left(
0,...,0,1,-1,...,1,0,...,0\right) ^{T}
\end{array}
\end{equation*}
Note that there exist $z_{2}$ pairs of consecutive points $y^{\left(
k_{w}\right) },\,y^{\left( k_{w}+1\right) \text{ }}$between $y^{\left(
i_{2}\right) }$ and $y^{\left( r\right) }$, with $1\leq w\leq z_{2}$.
Similarly and by definition of $U$ we verify that 
\begin{equation*}
U\,V_{\beta }\,\left( y^{\left( i_{2}\right) }\right) ^{T}=\left(
0,...,0,\varepsilon \left( y^{\left( i_{2}\right) }\right) \,\left|
F^{\prime }\left( y^{\left( i_{2}\right) }\right) \right| ^{-\beta
},0,...,0\right) ^{T}
\end{equation*}
where $\varepsilon \left( y^{\left( i_{2}\right) }\right) \,\left| F^{\prime
}\left( y^{\left( i_{2}\right) }\right) \right| ^{-\beta }$ is in the $j^{%
\text{th}}$-position. As the weight is constant on each interval $I_{j}$,
the desired result follows. Nevertheless, according to the behavior of $F$
several different cases may occur. The proof of the remaining cases is very
similar to the previous one.
\end{proof}

The above lemma suggest the next definition and result. Associated to each
matrix $V_{\beta }$ we have only one map $\xi :H_{0}\rightarrow H_{0}$,
which reflects the monotonicity of $F$. The map $\xi $ is represented by the 
$n\times n$ weighted diagonal matrix $K_{\beta }=\left[ k_{ij}\right] $,
where 
\begin{equation*}
k_{ii}:=\varepsilon \left( I_{i}\right) \,\left| F^{\prime }\left(
I_{i}\right) \right| ^{-\beta }
\end{equation*}
with $\varepsilon \left( I_{i}\right) =sign\left( F^{\prime }\left( x\right)
\right) $, $x\in I_{i}$.

\begin{lemma}
The next diagram is commutative 
\begin{equation*}
\begin{array}{ccccc}
&  & U &  &  \\ 
& C_{0} & \longrightarrow & H_{0} &  \\ 
V_{\beta } & \downarrow &  & \downarrow & K_{\beta } \\ 
& C_{0} & \longrightarrow & H_{0} &  \\ 
&  & U &  & 
\end{array}
.
\end{equation*}
\end{lemma}

The main results can now be stated. The next theorems establish the relation
between the weighted transition matrix, the weighted matrix $V_{\beta }$ and
the weighted kneading determinant.

\begin{theorem}
Under the conditions of the previous lemmas, we have the following relation
between the characteristic polynomials of the matrices $Q_{\beta
},\,V_{\beta }$ and $K_{\beta }$%
\begin{equation*}
P_{V_{\beta }}\left( t\right) =\det \left( I-t\,V_{\beta }\right)
=P_{Q_{\beta }}\left( t\right) \,P_{K_{\beta }}\left( t\right) .
\end{equation*}
\end{theorem}

\begin{proof}
The statement is a consequence of the above lemmas and according to some
homological algebra results, \cite{L}.
\end{proof}

We define a permutation matrix $\Pi $ associated to the permutation $\rho $
above defined, that maps the system of vectors $\left( z^{\left( 1\right)
},\,z^{\left( 2\right) },...,\,z^{\left( q\right) }\right) $ in to the
system of vectors $\left( y^{\left( 1\right) },\,y^{\left( 2\right)
},...,\,y^{\left( q\right) }\right) $. Using the weighted matrix $V_{\beta }$
and the permutation matrix $\Pi $, we define a new weighted matrix $\Theta
_{\beta }$ through the next equality 
\begin{equation*}
\Theta _{\beta }:=\Pi \,V_{\beta }\,\Pi ^{T}.
\end{equation*}
The matrix $\Pi $ is invertible, consequently we have $P_{\Theta _{\beta
}}\left( t\right) =P_{V_{\beta }}\left( t\right) $.

\begin{theorem}
If the kneading data associated to an expanding discontinuous map with holes 
$F$ corresponds to periodic, eventually periodic orbits or to orbits that
lie in the hole then the weighted kneading determinant is given by 
\begin{equation*}
D\left( t,\beta \right) =\frac{P_{Q_{\beta }}\left( t\right) }{R\left(
t\right) }
\end{equation*}
where $R\left( t\right) $ is a product of weighted cyclotomic polynomials
correspondent to those periodic or eventually periodic orbits.
\end{theorem}

It is obvious that this statement strongly depends on the number of laps and
the kneading data associated to $F$. For this reason, the analysis of the
general situation is difficult. We will prove the statement for a map $%
F=\left( F_{1},F_{2}\right) $. The general case follows in a similar way.

\begin{proof}
Consider $F=\left( F_{1},\,F_{2}\right) $, $\mathcal{P}_{I}=\left\{
I_{1},\,I_{2},\,I_{3}\right\} $ where $I_{2}$ is the hole, $F_{1}\left(
a_{2}\right) =1,\,F_{2}\left( a_{3}\right) =0$ and $\mathcal{A}=\left\{
L,\,H,\,R\right\} $. The orbits of the points $a_{1}^{+}$ and $a_{4}^{-}$
can be periodic, eventually periodic or lie in the hole. We consider that
the kneading data associated to this map is given by 
\begin{equation*}
\left( o\left( x^{\left( 2\right) }\right) ,\,o\left( x^{\left( 5\right)
}\right) \right) =\left( \left( L\,S_{1}^{\left( 2\right)
}...S_{p-1}^{\left( 2\right) }\right) ^{\infty },\,\left( M\,S_{1}^{\left(
5\right) }...S_{q-1}^{\left( 5\right) }\right) ^{\infty }\right) ,
\end{equation*}
where $p$ and $q$ are the periods of the orbits. The weighted kneading
increments are 
\begin{equation*}
\nu _{a_{2}}\left( t,\beta \right) =\theta _{a_{2}^{-}}\left( t,\beta
\right) \text{ and }\nu _{a_{3}}\left( t,\beta \right) =\theta
_{a_{3}^{+}}\left( t,\beta \right)
\end{equation*}
where 
\begin{equation*}
\theta _{a_{2}^{-}}\left( t,\beta \right) =\frac{L+\sum\limits_{i=1}^{p-1}%
\tau _{i}\left( x^{\left( 2\right) }\right) \text{\thinspace }%
t^{i}\,S_{i}^{\left( 2\right) }}{1-\tau _{p}\left( x^{\left( 2\right)
}\right) \,t^{p}}
\end{equation*}
and 
\begin{equation*}
\theta _{a_{3}^{+}}\left( t,\beta \right) =\frac{R+\sum\limits_{i=1}^{q-1}%
\tau _{i}\left( x^{\left( 5\right) }\right) \text{\thinspace }%
t^{i}\,S_{i}^{\left( 5\right) }}{1-\tau _{q}\left( x^{\left( 5\right)
}\right) \,t^{q}}.
\end{equation*}
If we write 
\begin{equation*}
L_{p}=\sum\limits_{\substack{ i=1  \\ \,S_{i}^{\left( 2\right) }=\,L}}%
^{p-1}\tau _{i}\left( x^{\left( 2\right) }\right) \text{\thinspace }t^{i}%
\text{ and }R_{p}=\sum\limits_{\substack{ i=1  \\ \,S_{i}^{\left( 2\right)
}=\,R}}^{p-1}\tau _{i}\left( x^{\left( 2\right) }\right) \text{\thinspace }%
t^{i}
\end{equation*}
and analogously for $L_{q}$ and $R_{q}$, then we have 
\begin{equation*}
\nu _{a_{2}}\left( t,\beta \right) =\frac{1+L_{p}}{1-\tau _{p}\left(
x^{\left( 2\right) }\right) \,t^{p}}\,L+\frac{R_{p}}{1-\tau _{p}\left(
x^{\left( 2\right) }\right) \,t^{p}}\,R
\end{equation*}
and 
\begin{equation*}
\nu _{a_{3}}\left( t,\beta \right) =\frac{L_{q}}{1-\tau _{q}\left( x^{\left(
5\right) }\right) \,t^{q}}\,L+\frac{1+R_{q}}{1-\tau _{q}\left( x^{\left(
5\right) }\right) \,t^{q}}\,R.
\end{equation*}
The weighted kneading determinant $D\left( t,\beta \right) $ for these
kneading data is 
\begin{equation*}
D\left( t,\beta \right) =\frac{1}{\left( 1-\tau _{p}\left( x^{\left(
2\right) }\right) \,t^{p}\right) \,\left( 1-\tau _{q}\left( x^{\left(
5\right) }\right) \,t^{q}\right) }\left| 
\begin{array}{cc}
1+L_{p} & R_{p} \\ 
&  \\ 
L_{q} & 1+R_{q}
\end{array}
\right| .
\end{equation*}

Let $\det \overline{\Theta }_{\beta }=\det (I-t\,\Theta _{\beta })$ be the
characteristic polynomial of the matrix $\Theta _{\beta }$, where $I$ is the
identity matrix. Thus, we have 
\begin{equation*}
\overline{\Theta }_{\beta }=\left[ 
\begin{array}{ccccc}
1-\mu _{2,0}t & -\mu _{2,1}t & \ldots & -\mu _{2,p-2}t & -\mu
_{2,p-1}t-\theta _{2,p-1}t \\ 
-\theta _{2,0}t & 1 & \ldots & 0 & 0 \\ 
\vdots & \vdots & \ddots & \vdots & \vdots \\ 
0 & 0 & \ldots & 1 & 0 \\ 
0 & 0 & \ldots & -\theta _{2,p-2}t & 1 \\ 
-\delta _{2,0}t & -\delta _{2,1}t & \ldots & -\delta _{2,p-2}t & -\delta
_{2,p-1}t \\ 
0 & 0 & \ldots & 0 & 0 \\ 
\vdots & \vdots & \ddots & \vdots & \vdots \\ 
0 & 0 & \ldots & 0 & 0 \\ 
0 & 0 & \ldots & 0 & 0
\end{array}
\right.
\end{equation*}
\begin{equation*}
\left. 
\begin{array}{ccccc}
-\mu _{5,0}t & -\mu _{5,1}t & \ldots & -\mu _{5,q-2}t & -\mu _{5,q-1}t \\ 
0 & 0 & \ldots & 0 & 0 \\ 
\vdots & \vdots & \ddots & \vdots & \vdots \\ 
0 & 0 & \ldots & 0 & 0 \\ 
0 & 0 & \ldots & 0 & 0 \\ 
1-\delta _{5,0}t & -\delta _{5,1}t & \ldots & -\delta _{5,q-2}t & -\delta
_{5,q-1}t-\theta _{5,q-1}t \\ 
-\theta _{5,0}t & 1 & \ldots & 0 & 0 \\ 
\vdots & \vdots & \ddots & \vdots & \vdots \\ 
0 & 0 & \ldots & 1 & 0 \\ 
0 & 0 & \ldots & -\theta _{5,q-2}t & 1
\end{array}
\right] .
\end{equation*}
where 
\begin{eqnarray*}
\theta _{2,k_{1}} &=&\varepsilon \left( x_{k_{1}}^{\left( 2\right) }\right)
\left| F^{\prime }\left( x_{k_{1}}^{\left( 2\right) }\right) \right|
^{-\beta } \\
\mu _{2,k_{1}} &\in &\left\{ 0,\,\pm \varepsilon \left( x_{k_{1}}^{\left(
2\right) }\right) \left| F^{\prime }\left( x_{k_{1}}^{\left( 2\right)
}\right) \right| ^{-\beta }\right\} \text{ with }0\leq k_{1}\leq p-1,
\end{eqnarray*}
\begin{eqnarray*}
\theta _{5,k_{2}} &=&\varepsilon \left( x_{k_{2}}^{\left( 5\right) }\right)
\left| F^{\prime }\left( x_{k_{2}}^{\left( 5\right) }\right) \right|
^{-\beta } \\
\mu _{5,k_{2}} &\in &\left\{ 0,\,\pm \varepsilon \left( x_{k_{2}}^{\left(
5\right) }\right) \left| F^{\prime }\left( x_{k_{2}}^{\left( 5\right)
}\right) \right| ^{-\beta }\right\} \text{ with }0\leq k_{2}\leq q-1
\end{eqnarray*}
and similarly for $\delta _{2,k_{1}}$ and $\delta _{5,k_{2}}$. Using matrix
elementary operations for the matrix $I-t\,\Theta _{\beta }$, we have the
following equivalent matrix 
\begin{equation*}
\left[ 
\begin{array}{cc}
1-\sum\limits_{k=0}^{p-1}\mu _{2,k}\tau _{k}\left( x_{0}^{\left( 2\right)
}\right) t^{k+1}-\tau _{p}\left( x_{0}^{\left( 2\right) }\right) t^{p} & 
-\sum\limits_{k=0}^{q-1}\mu _{5,k}\tau _{k}\left( x_{0}^{\left( 5\right)
}\right) t^{k+1} \\ 
&  \\ 
-\sum\limits_{k=0}^{p-1}\delta _{2,k}\tau _{k}\left( x_{0}^{\left( 2\right)
}\right) t^{k+1} & 1-\sum\limits_{k=0}^{q-1}\delta _{5,k}\tau _{k}\left(
x_{0}^{\left( 5\right) }\right) t^{k+1}-\tau _{q}\left( x_{0}^{\left(
5\right) }\right) t^{q}
\end{array}
\right] .
\end{equation*}
Now we will compare the elements of the above matrix to the elements of the
weighted kneading matrix $N\left( t,\beta \right) .$ Note that 
\begin{equation}
\begin{array}{l}
1-\sum\limits_{k=0}^{p-1}\mu _{2,k}\,\tau _{k}\left( x_{0}^{\left( 2\right)
}\right) \,t^{k+1}-\tau _{p}\left( x_{0}^{\left( 2\right) }\right) \,t^{p}
\\ 
\\ 
=\left( 1-\mu _{2,0}\,t\right) +\sum\limits_{k=1}^{p-1}\tau _{k}\left(
x_{0}^{\left( 2\right) }\right) t^{k}\left( -\mu _{2,k}\,t\right) \\ 
\\ 
+\tau _{p-1}\left( x_{0}^{\left( 2\right) }\right) t^{p-1}\left(
-\varepsilon \left( x_{p-1}^{\left( 2\right) }\right) \left| F^{\prime
}\left( x_{p-1}^{\left( 2\right) }\right) \right| ^{-\beta }t\right) .
\end{array}
\label{FK2.1}
\end{equation}
On the other hand 
\begin{equation}
\begin{array}{l}
\left( 1-\varepsilon \left( I_{1}\right) \,\left| F^{\prime }\left(
I_{1}\right) \right| ^{-\beta }t\right) \left( 1+L_{p}\right) \\ 
\\ 
=\left( 1-\varepsilon \left( I_{1}\right) \,\left| F^{\prime }\left(
I_{1}\right) \right| ^{-\beta }t\right) +\sum\limits_{\substack{ k=1  \\ %
\,S_{k}^{\left( 2\right) }=\,L}}^{p-1}\tau _{k}\left( x_{0}^{\left( 2\right)
}\right) \text{\thinspace }t^{k} \\ 
\\ 
+\sum\limits_{\substack{ k=1  \\ \,S_{k}^{\left( 2\right) }=\,L}}^{p-1}\tau
_{k}\left( x_{0}^{\left( 2\right) }\right) \text{\thinspace }t^{k}\left(
-\varepsilon \left( I_{1}\right) \,\left| F^{\prime }\left( I_{1}\right)
\right| ^{-\beta }t\right) .
\end{array}
\label{FK2.2}
\end{equation}
If $\mu _{2,k}\neq 0$ with $1\leq k\leq p-1$ then $\tau _{k+1}\left(
x_{0}^{\left( 2\right) }\right) \neq 0$ on $L_{p}$. That is, 
\begin{equation*}
\tau _{k}\left( x_{0}^{\left( 2\right) }\right) \text{\thinspace }%
t^{k}\left( \pm \varepsilon \left( x_{k}^{\left( 2\right) }\right) \left|
F^{\prime }\left( x_{k}^{\left( 2\right) }\right) \right| ^{-\beta }t\right)
=\pm \tau _{k+1}\left( x_{0}^{\left( 2\right) }\right) \text{\thinspace }%
t^{k+1}.
\end{equation*}
In particular, if $\mu _{2,p-1}\neq 0$ then $x_{p-1}^{\left( 2\right) }$ is
associated to the symbol $R$, i.e., 
\begin{equation*}
\mu _{2,p-1}=-\varepsilon \left( x_{p-1}^{\left( 2\right) }\right) \left|
F^{\prime }\left( x_{p-1}^{\left( 2\right) }\right) \right| ^{-\beta }.
\end{equation*}
Consequently, in (\ref{FK2.1}) we have 
\begin{equation*}
\tau _{p-1}\left( x_{0}^{\left( 2\right) }\right) \text{\thinspace }%
t^{p-1}\varepsilon \left( x_{p-1}^{\left( 2\right) }\right) \left| F^{\prime
}\left( x_{p-1}^{\left( 2\right) }\right) \right| ^{-\beta }t=\tau
_{p}\left( x_{0}^{\left( 2\right) }\right) \text{\thinspace }t^{p}.
\end{equation*}
Hence, in (\ref{FK2.2}), the fact that the orbit is periodic, implies that
we return to the symbol $L$. Thus, we have 
\begin{equation*}
1-\mu _{2,0}\,t=1-\varepsilon \left( I_{1}\right) \,\left| F^{\prime }\left(
I_{1}\right) \right| ^{-\beta }t.
\end{equation*}

Let $R\left( t\right) $ be the product of cyclotomic polynomials and $%
P_{K_{\beta }}\left( t\right) $ be the characteristic polynomial of the
matrix $K_{\beta }$ associated to $F$. Set 
\begin{equation*}
D^{\ast }\left( t,\beta \right) =\left| 
\begin{array}{cc}
\left( 1-\varepsilon \left( I_{1}\right) \,\left| F^{\prime }\left(
I_{1}\right) \right| ^{-\beta }t\right) \left( 1+L_{p}\right) & \left(
1-\varepsilon \left( I_{2}\right) \,\left| F^{\prime }\left( I_{2}\right)
\right| ^{-\beta }t\right) R_{p} \\ 
&  \\ 
\left( 1-\varepsilon \left( I_{1}\right) \,\left| F^{\prime }\left(
I_{1}\right) \right| ^{-\beta }t\right) L_{q} & \left( 1-\varepsilon \left(
I_{2}\right) \,\left| F^{\prime }\left( I_{2}\right) \right| ^{-\beta
}t\right) \left( 1+R_{q}\right)
\end{array}
\right| .
\end{equation*}
Using the above comparison between the elements of the equivalent matrix to $%
I-t\,\Theta _{\beta }$ and the elements of the weighted kneading matrix, we
have 
\begin{equation*}
D\left( t,\beta \right) =\frac{D^{\ast }\left( t,\beta \right) }{R\left(
t\right) \,P_{K_{\beta }}\left( t\right) }=\det \left( I-t\,\Theta _{\beta
}^{T}\right) =P_{V_{\beta }}\left( t\right) .
\end{equation*}
By Theorem 1, the desired result follows.
\end{proof}

The following statement will allows us to compute explicitly the Hausdorff
dimension, the escape rate and the topological entropy.

\begin{theorem}
Let $D\left( t,\beta \right) $ be the weighted kneading determinant, under
the conditions of the previous theorem.

(i) If $\beta $ is the unique solution of $D\left( 1,\beta \right) =0$, then 
$\beta $ is the Hausdorff dimension of the attractor $E$.

(ii) If $t_{1}$ is the least real positive solution of $D\left( t,1\right)
=0 $, then $\log \left( t_{1}\right) $ is the escape rate of the pair $%
\left( E,F\right) $.

(iii) If $t_{0}$ is the least real positive solution of $D\left( t,0\right)
=0$, then $\log \left( t_{0}^{-1}\right) $ is the topological entropy of the
map $F$.
\end{theorem}

\begin{proof}
Considering the transfer operator given in (\ref{FK1.2}), we have 
\begin{equation*}
\left( L_{\phi _{j}}\,g\right) \left( x\right) =\sum\limits_{j=1}^{m}\left|
F_{j}^{\prime }\left( x\right) \right| ^{-\beta }\,g\left( f_{j}\left(
x\right) \right) \,\chi _{F\left( int\,J_{j}\right) }.
\end{equation*}
Let $a_{ij}$ be the entries of the transition matrix $A$. For each $J_{i}\in 
\mathcal{P}_{I}^{\prime }$, with $1\leq i\leq m$, and $\beta \in \mathbb{R}$
the eigenvalue equation corresponding to an eigenvalue $\lambda _{\beta }$
is 
\begin{equation*}
\sum\limits_{j=1}^{m}\frac{a_{ij}}{\left| F_{j}^{\prime }\left( x\right)
\right| ^{\beta }}\,v_{j}=\lambda _{\beta }\,v_{i}
\end{equation*}
for the operator $L_{\phi }$ characterized by the matrix $Q_{\beta }.$
According to \cite{R1} and using (\ref{FK1.4}), the largest eigenvalue of
the transfer operator is $\exp P\left( \beta \right) $. Hence, $\exp P\left(
\beta \right) $ is the spectral radius $\lambda _{\beta }$ of the matrix $%
Q_{\beta }$.

If $\beta $ is the unique solution of $D\left( 1,\beta \right) =0$, then by
Theorem 2 and (\ref{FK1.4}), we get $P\left( \beta \right) =0$. By \cite{F}
and \cite{LS2}, we can conclude that $\beta =\dim _{H}\left( E\right) $.

On the other hand, considering the parameter $\beta =1$, we have that $%
\lambda _{1}=\exp P\left( 1\right) $ is the largest eigenvalue of the matrix 
$Q_{1}$. The second statement follows from \cite{BC}, where the escape rate $%
\gamma $ is given by $\gamma =-P\left( 1\right) $. Thus, the escape rate is $%
\gamma =\log \left( \lambda _{1}^{-1}\right) $, where $\lambda
_{1}^{-1}=t_{1}$ is the least real positive solution of $P_{Q_{1}}\left(
t\right) =0$.

If $\beta =0$, then the determinant $D\left( t,0\right) $ corresponds to the
kneading determinant described in \cite{LS1}, where $t_{0}^{-1}=\lambda _{0}$
is the growth number of $F$, i.e., the spectral radius of the transition
matrix $A$. Consequently, $\log \left( \lambda _{0}\right) $ is the
topological entropy of the map $F$.
\end{proof}

\begin{remark}
The theory presented in this paper with respect to periodic, eventually
periodic orbits or to the orbits that lie in the holes is also valid for
aperiodic orbits. In this case, the invariant coordinates associated to the
turning points and to the discontinuity points are formal power series. The
computation of the topological invariants is done by approximation using:
periodic, eventually periodic orbits or the orbits that lie in the holes.
\end{remark}

The above results are illustrated in the next example, showing in detail the
techniques under discussion.

\begin{example}
Let 
\begin{equation*}
F\left( x\right) =\left\{ 
\begin{array}{ll}
\frac{x}{a} & \text{if }x\in \left[ 0,\frac{1}{8}\right] \\ 
&  \\ 
-\frac{x}{a}+1 & \text{if }x\in \left[ \frac{1}{8},\frac{1}{4}\right] \\ 
&  \\ 
\frac{x}{b}-\frac{3}{4} & \text{if }x\in \left[ \frac{1}{4},\frac{7}{12}%
\right] \\ 
&  \\ 
\frac{x}{c}-1 & \text{if }x\in \left[ \frac{4}{6},1\right]
\end{array}
\right.
\end{equation*}
with $a=1/4,\,b=1/3$ and $c=1/2.$ Considering the orbits of the points $%
a_{2}^{\pm },\,a_{3}^{\pm },\,a_{4}^{-}$ and $a_{5}^{+}$, the kneading data
to the map $F$ are 
\begin{equation*}
\left( L\left( M_{2}R\right) ^{\infty },\,M_{1}\left( M_{2}R\right) ^{\infty
},\,M_{1}L^{\infty },\,M_{2}L^{\infty },\,M_{2}R^{\infty
},\,RM_{2}M_{2}L^{\infty }\right) .
\end{equation*}
The weighted invariant coordinates of each point are 
\begin{equation*}
\theta _{a_{2}^{-}}\left( t,\beta \right) =L+\frac{a^{\beta }\,t}{1-\left(
b\,c\right) ^{\beta }\,t^{2}}\,M_{2}+\frac{\left( a\,b\right) ^{\beta
}\,t^{2}}{1-\left( b\,c\right) ^{\beta }\,t^{2}}\,R
\end{equation*}
\begin{equation*}
\theta _{a_{2}^{+}}\left( t,\beta \right) =M_{1}-\frac{a^{\beta }\,t}{%
1-\left( b\,c\right) ^{\beta }\,t^{2}}\,M_{2}-\frac{\left( a\,b\right)
^{\beta }\,t^{2}}{1-\left( b\,c\right) ^{\beta }\,t^{2}}\,R
\end{equation*}
\begin{equation*}
\theta _{a_{3}^{-}}\left( t,\beta \right) =-\frac{a^{\beta }\,t}{1-a^{\beta
}\,t}\,L+M_{1};\,\theta _{a_{3}^{+}}\left( t,\beta \right) =\frac{b^{\beta
}\,t}{1-a^{\beta }\,t}\,L+M_{2}
\end{equation*}
\begin{equation*}
\theta _{a_{4}^{-}}\left( t,\beta \right) =M_{2}+\frac{b^{\beta }\,t}{%
1-c^{\beta }\,t}\,R;\,\theta _{a_{5}^{+}}\left( t,\beta \right) =\frac{%
(b^{2}c)^{\beta }\,t^{3}}{1-a^{\beta }\,t}\,L+\left( c^{\beta }\,t+\left(
b\,c\right) ^{\beta }\,t^{2}\right) \,M_{2}+R.
\end{equation*}
Consequently, the weighted kneading determinant is 
\begin{equation*}
\begin{array}{lll}
D\left( t,\beta \right) & = & \left| 
\begin{array}{cccc}
-1 & 1 & \frac{-2a^{\beta }\,t}{1-\left( b\,c\right) ^{\beta }\,t^{2}} & 
\frac{-2\left( a\,b\right) ^{\beta }\,t^{2}}{1-\left( b\,c\right) ^{\beta
}\,t^{2}} \\ 
\frac{a^{\beta }\,t+b^{\beta }\,t}{1-a^{\beta }\,t} & -1 & 1 & 0 \\ 
0 & 0 & 1 & \frac{b^{\beta }\,t}{1-c^{\beta }\,t} \\ 
\frac{(b^{2}c)^{\beta }\,t^{3}}{1-a^{\beta }\,t} & 0 & c^{\beta }\,t+\left(
b\,c\right) ^{\beta }\,t^{2} & 1
\end{array}
\right| \\ 
&  &  \\ 
& = & \frac{1}{\left( 1-a^{\beta }\,t\right) \,\left( 1-c^{\beta }\,t\right)
\,\left( 1-\left( b\,c\right) ^{\beta }\,t^{2}\right) }\left( 1-(2a^{\beta
}+b^{\beta }+c^{\beta })\,t\right. \\ 
&  &  \\ 
&  & +(2\left( a\,c\right) ^{\beta }-\left( b\,c\right) ^{\beta
})\,t^{2}+\left( 4\left( a\,b\,c\right) ^{\beta }+\left( b^{2}c\right)
^{\beta }+\left( b\,c^{2}\right) ^{\beta }\right) \,t^{3} \\ 
&  &  \\ 
&  & \left. +\left( 2\left( a\,b^{2}c\right) ^{\beta }-2\left(
a\,b\,c^{2}\right) ^{\beta }\right) \,t^{4}-2\left( a\,b^{2}c^{2}\right)
^{\beta }\,t^{5}\right) .
\end{array}
\end{equation*}

The orbits of the points $a_{2}^{\pm },\,a_{3}^{\pm },\,a_{4}^{-}$ and $%
a_{5}^{+}$ determine a Markov partition of $\left[ 0,1\right] $, $\mathcal{P}%
_{I}^{\prime }=\left\{ I_{1},...,\,I_{8}\right\} $ where $I_{6}$ is the
hole. The matrices correspondent to this map are 
\begin{equation*}
B=\left[ 
\begin{array}{rrrrrrrr}
-1 & 0 & 0 & 0 & 0 & 0 & 0 & 0 \\ 
1 & 0 & 0 & 0 & 0 & 0 & 0 & 0 \\ 
0 & -1 & 0 & 0 & 0 & 0 & 0 & 0 \\ 
0 & 1 & 0 & 0 & 0 & 0 & 0 & 0 \\ 
0 & 0 & -1 & 0 & 0 & 0 & 0 & 0 \\ 
0 & 0 & 1 & -1 & 0 & 0 & 0 & 0 \\ 
0 & 0 & 0 & 1 & -1 & 0 & 0 & 0 \\ 
0 & 0 & 0 & 0 & 1 & 0 & 0 & 0 \\ 
0 & 0 & 0 & 0 & 0 & -1 & 0 & 0 \\ 
0 & 0 & 0 & 0 & 0 & 1 & 0 & 0 \\ 
0 & 0 & 0 & 0 & 0 & 0 & -1 & 0 \\ 
0 & 0 & 0 & 0 & 0 & 0 & 1 & -1 \\ 
0 & 0 & 0 & 0 & 0 & 0 & 0 & 1
\end{array}
\right]
\end{equation*}
\begin{equation*}
Q_{\beta }=\left[ 
\begin{array}{cccccccc}
a^{\beta } & a^{\beta } & b^{\beta } & 0 & 0 & 0 & 0 & 0 \\ 
a^{\beta } & a^{\beta } & b^{\beta } & 0 & 0 & 0 & 0 & 0 \\ 
a^{\beta } & a^{\beta } & 0 & b^{\beta } & 0 & 0 & 0 & 0 \\ 
a^{\beta } & a^{\beta } & 0 & b^{\beta } & 0 & 0 & c^{\beta } & 0 \\ 
0 & 0 & 0 & b^{\beta } & 0 & 0 & 0 & c^{\beta } \\ 
0 & 0 & 0 & b^{\beta } & 0 & 0 & 0 & c^{\beta } \\ 
0 & 0 & 0 & b^{\beta } & 0 & 0 & 0 & c^{\beta } \\ 
0 & 0 & 0 & 0 & b^{\beta } & 0 & 0 & c^{\beta }
\end{array}
\right]
\end{equation*}
\begin{equation*}
V_{\beta }=\left[ 
\begin{array}{ccccccccccccc}
a^{\beta } & 0 & 0 & -a^{\beta } & b^{\beta } & 0 & 0 & 0 & 0 & 0 & 0 & 0 & 0
\\ 
0 & a^{\beta } & 0 & a^{\beta } & -b^{\beta } & 0 & 0 & 0 & 0 & 0 & 0 & 0 & 0
\\ 
0 & -a^{\beta } & 0 & -a^{\beta } & b^{\beta } & 0 & 0 & 0 & 0 & 0 & 0 & 0 & 
0 \\ 
0 & a^{\beta } & -a^{\beta } & 0 & -b^{\beta } & 0 & 0 & 0 & 0 & 0 & 0 & 0 & 
0 \\ 
0 & -a^{\beta } & a^{\beta } & 0 & b^{\beta } & b^{\beta } & 0 & 0 & 0 & 0 & 
0 & 0 & 0 \\ 
0 & 0 & 0 & 0 & 0 & 0 & 0 & 0 & 0 & 0 & c^{\beta } & 0 & 0 \\ 
0 & a^{\beta } & -a^{\beta } & 0 & 0 & 0 & 0 & 0 & 0 & 0 & 0 & c^{\beta } & 0
\\ 
0 & 0 & 0 & 0 & 0 & 0 & b^{\beta } & b^{\beta } & 0 & 0 & -c^{\beta } & 
-c^{\beta } & 0 \\ 
0 & 0 & 0 & 0 & 0 & 0 & -b^{\beta } & -b^{\beta } & 0 & 0 & c^{\beta } & 
c^{\beta } & 0 \\ 
0 & 0 & 0 & 0 & 0 & 0 & b^{\beta } & b^{\beta } & 0 & 0 & -c^{\beta } & 
-c^{\beta } & 0 \\ 
0 & 0 & 0 & 0 & 0 & 0 & -b^{\beta } & -b^{\beta } & 0 & 0 & c^{\beta } & 
c^{\beta } & 0 \\ 
0 & 0 & 0 & 0 & 0 & 0 & b^{\beta } & 0 & 0 & 0 & 0 & 0 & 0 \\ 
0 & 0 & 0 & 0 & 0 & 0 & 0 & b^{\beta } & 0 & 0 & 0 & 0 & c^{\beta }
\end{array}
\right]
\end{equation*}
\begin{equation*}
U=\left[ 
\begin{array}{ccccccccccccc}
1 & 1 & 0 & 0 & 0 & 0 & 0 & 0 & 0 & 0 & 0 & 0 & 0 \\ 
0 & 0 & 1 & 1 & 0 & 0 & 0 & 0 & 0 & 0 & 0 & 0 & 0 \\ 
0 & 0 & 0 & 0 & 1 & 1 & 1 & 1 & 0 & 0 & 0 & 0 & 0 \\ 
0 & 0 & 0 & 0 & 0 & 0 & 0 & 0 & 0 & 0 & 1 & 1 & 1
\end{array}
\right] \,K_{\beta }=\left[ 
\begin{array}{cccc}
a^{\beta } & 0 & 0 & 0 \\ 
0 & -a^{\beta } & 0 & 0 \\ 
0 & 0 & b^{\beta } & 0 \\ 
0 & 0 & 0 & c^{\beta }
\end{array}
\right]
\end{equation*}
\begin{equation*}
\Pi =\left[ 
\begin{array}{ccccccccccccc}
1 & 2 & 3 & 4 & 5 & 6 & 7 & 8 & 9 & 10 & 11 & 12 & 13 \\ 
6 & 1 & 4 & 5 & 7 & 11 & 2 & 8 & 12 & 13 & 10 & 3 & 9
\end{array}
\right]
\end{equation*}
\begin{equation*}
\Theta _{\beta }=\left[ 
\begin{array}{ccccccccccccc}
a^{\beta } & 0 & 0 & 0 & a^{\beta } & 0 & -b^{\beta } & 0 & 0 & 0 & 0 & 0 & 0
\\ 
a^{\beta } & 0 & c^{\beta } & -a^{\beta } & 0 & 0 & 0 & 0 & 0 & 0 & 0 & 0 & 0
\\ 
0 & b^{\beta } & 0 & 0 & 0 & 0 & 0 & 0 & 0 & 0 & 0 & 0 & 0 \\ 
-a^{\beta } & 0 & 0 & 0 & -a^{\beta } & 0 & b^{\beta } & 0 & 0 & 0 & 0 & 0 & 
0 \\ 
a^{\beta } & 0 & 0 & -a^{\beta } & 0 & 0 & -b^{\beta } & 0 & 0 & 0 & 0 & 0 & 
0 \\ 
0 & 0 & 0 & 0 & -a^{\beta } & a^{\beta } & b^{\beta } & 0 & 0 & 0 & 0 & 0 & 0
\\ 
-a^{\beta } & 0 & 0 & a^{\beta } & 0 & 0 & b^{\beta } & 0 & 0 & 0 & b^{\beta
} & 0 & 0 \\ 
0 & b^{\beta } & -c^{\beta } & 0 & 0 & 0 & 0 & b^{\beta } & 0 & -c^{\beta }
& 0 & 0 & 0 \\ 
0 & 0 & 0 & 0 & 0 & 0 & 0 & b^{\beta } & c^{\beta } & 0 & 0 & 0 & 0 \\ 
0 & -b^{\beta } & c^{\beta } & 0 & 0 & 0 & 0 & -b^{\beta } & 0 & c^{\beta }
& 0 & 0 & 0 \\ 
0 & 0 & 0 & 0 & 0 & 0 & 0 & 0 & 0 & c^{\beta } & 0 & 0 & 0 \\ 
0 & -b^{\beta } & c^{\beta } & 0 & 0 & 0 & 0 & -b^{\beta } & 0 & c^{\beta }
& 0 & 0 & 0 \\ 
0 & b^{\beta } & -c^{\beta } & 0 & 0 & 0 & 0 & b^{\beta } & 0 & -c^{\beta }
& 0 & 0 & 0
\end{array}
\right] .
\end{equation*}
The relation between the characteristic polynomials of the matrices $%
V_{\beta },\,Q_{\beta }$ and $K_{\beta }$ is 
\begin{equation*}
\begin{array}{lll}
P_{V_{\beta }}\left( t\right) & = & \left( 1-(2a^{\beta }+b^{\beta
}+c^{\beta })t+(2\left( a\,c\right) ^{\beta }-\left( b\,c\right) ^{\beta
})t^{2}+\left( 4\left( a\,b\,c\right) ^{\beta }+\left( b^{2}c\right) ^{\beta
}\right) t^{3}\right. \\ 
&  &  \\ 
&  & \left. +\left( b\,c^{2}\right) ^{\beta }\,t^{3}+\left( 2\left(
a\,b^{2}c\right) ^{\beta }-2\left( a\,b\,c^{2}\right) ^{\beta }\right)
\,t^{4}-2\left( a\,b^{2}c^{2}\right) ^{\beta }\,t^{5}\right) \\ 
&  &  \\ 
&  & \left( 1-a^{\beta }\,t\right) \,\left( 1+a^{\beta }\,t\right) \,\left(
1-b^{\beta }\,t\right) \,\left( 1-c^{\beta }\,t\right) \\ 
&  &  \\ 
& = & P_{Q_{\beta }}\left( t\right) \,P_{K_{\beta }}\left( t\right) .
\end{array}
\end{equation*}
Then we can verify the statements of Theorem 1 and 2. By Theorem 3, we have 
\begin{equation*}
\dim _{H}\left( E\right) =0.91994...;\,\gamma =0.0877769...\text{ and }%
h_{top}=1.11531...\text{.}
\end{equation*}
\end{example}

\end{document}